# Tiling of Strip Lattices and Asymptotics

## Valcho Milchev

**Abstract.** This paper examines the tilings of a strip with equilateral triangles in a $2 \times n$ triangular strip. The number of ways in which the lattices can be covered with a combination of tiles of the two types of triangles is related to Pell numbers and Pell-Lucas numbers. There are multiple dependencies among the tilings of the lattices. Additionally, the question of the number of tiles required for all possible tilings – both the number of tiles in total and by type – is developed for the first time in this paper. Attention is given to asymptotics as well.

We consider an unlimited number of tiles in the shape of equilateral triangles of two sizes – triangles with a side length of 1 (referred to as small tiles or small triangles) and triangles with a side length of 2 (referred to as large tiles or large triangles). They are shown in Figure 1. It is clear that a large triangle consists of four small triangles. ([1]).

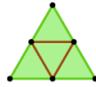

**Figure 1**

Let us examine the triangular lattices shown in Figure 2. They are strips composed of small equilateral triangles placed next to each other. We analyze the number of ways in which these lattices can be covered using a combination of the two types of triangular tiles.

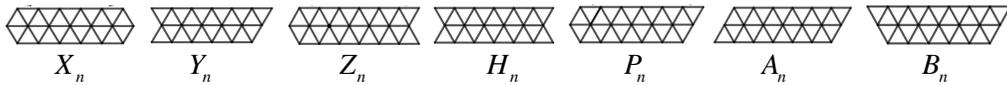

$X_n$    $Y_n$    $Z_n$    $H_n$    $P_n$    $A_n$    $B_n$

**Figure 2**

For convenience, the number of tilings of the lattices is denoted exactly the same as the lattices themselves: $H_n$, $P_n$, $X_n$, $Y_n$, $Z_n$, $A_n$, $B_n$.

**Theorem 1.** The number of ways in which the lattices $H_n$ and $P_n$ can be covered with a combination of the two types of triangular tiles corresponds to the Pell numbers and the Pell-Lucas numbers respectively (Lattice $P_n$ is interpreted as being obtained from lattice $H_n$ by removing one small triangle - Figure 3).

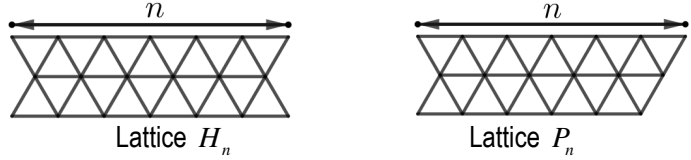

**Figure 3**

**Proof.** We derive recurrence relations for the number of tilings of the two strip lattices $H_n$ and $P_n$. It is easy to count that $H_1 = 1$, $H_2 = 3$, $P_1 = 1$, $P_2 = 2$ (Figure 4).

From the two obtained recurrence equations (Figure 5) $H_n = H_{n-1} + 2P_{n-1}$ and

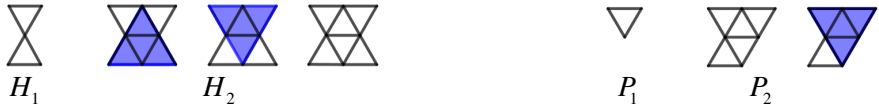

**Figure 4 :** $H_1 = 1$, $H_2 = 3$, $P_1 = 1$, $P_2 = 2$.

$P_n = H_{n-1} + P_{n-1}$ we find that for $n \geq 3$:
$$H_n = 2H_{n-1} + H_{n-2}, \quad P_n = 2P_{n-1} + P_{n-2}.$$

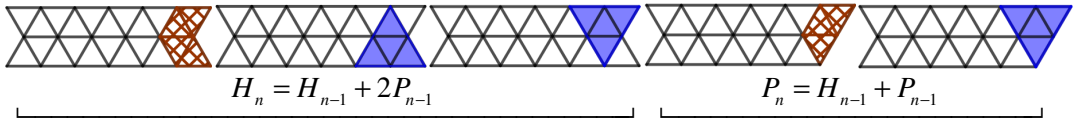

**Figure 5.** The equations $H_n = H_{n-1} + 2P_{n-1}$, $P_n = H_{n-1} + P_{n-1}$ are illustrated.

Thus $\{P_n\} = \{1, 2, 5, 12, 29, 70, ...\}$, $\{H_n\} = \{1, 3, 7, 17, 41, 99, ...\}$.

Now from the recurrent dependencies for $H_n$ and $P_n$ we can determine that $P_0 = 0$ and $H_0 = 1$.

It has become clear that $\{P_n\} = \{1, 2, 5, 12, 29, 70, ...\}$ are the famous Pell numbers ([2], **A000129**) whereas $\{H_n\} = \{1, 3, 7, 17, 41, 99, ...\}$ are a sequence of numbers ([2], **A001333**). included in The On-Line Encyclopedia of Integer Sequences.

These are the generating functions of the linear recurrence sequences $H_n$ and $P_n$:
$$G_H(x) = \frac{1-x}{1-2x-x^2}, \quad G_P(x) = \frac{x}{1-2x-x^2}. \square$$

The next statement is about tilings of a strip for the lattices $X_n$, $Y_n$ and $Z_n$. These triangular lattices are shown in Figure 6.

**Theorem 2.** The number of ways in which the lattices $X_n$, $Y_n$ and $Z_n$ can be covered with a combination of the two types of triangular tiles corresponds to the Pell numbers and the Pell-Lucas numbers respectively.

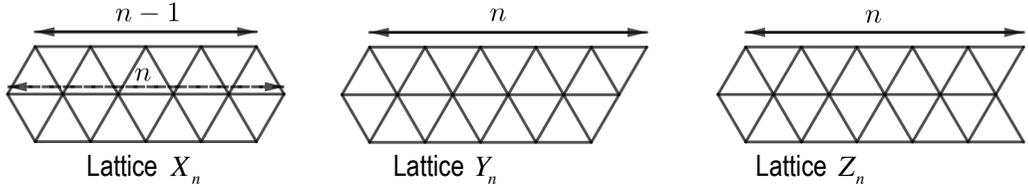

**Figure 6**

**Proof.** We derive recurrence relations for the number of tilings of the two strip lattices $Y_n$ and $Z_n$. It is easy to count that $Z_1 = 1$, $Z_2 = 3$, $Y_1 = 1$, $Y_2 = 2$ (Figure 7).

From the two obtained recurrence equations $Z_n = Z_{n-1} + 2Y_{n-1}$ and $Y_n = Z_{n-1} + Y_{n-1}$

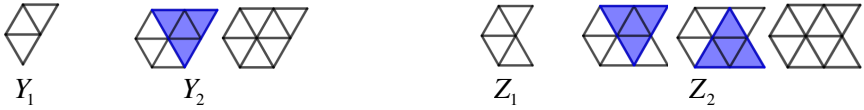

**Figure 7**

(Figure 8) we find that for $n \geq 3$: $Z_n = 2Z_{n-1} + Z_{n-2}$, $Y_n = 2Y_{n-1} + Y_{n-2}$. Thus $\{Z_n\} = \{1, 3, 7, 17, 41, 99, ...\}$ and $\{Y_n\} = \{1, 2, 5, 12, 29, 70, ...\}$.

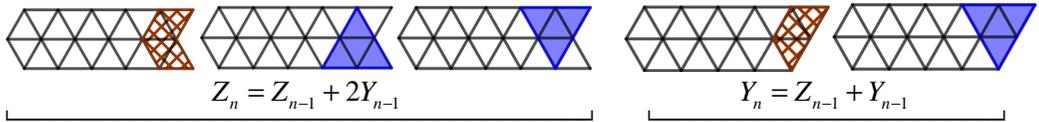

**Figure 8.** The equations $Z_n = Z_{n-1} + 2Y_{n-1}$, $Y_n = Z_{n-1} + Y_{n-1}$ are illustrated.

Let us examine the triangular lattices $X_n$. From the two obtained reccurence equations from Figure 8 and Figure 9 $Y_n = Z_{n-1} + Y_{n-1}$, $Y_n = Y_{n-1} + X_n$, we can see that $X_n = Z_{n-1}$.

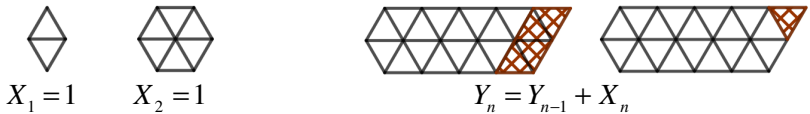

**Figure 9.**

From the two obtained reccurence equations we conclude that $\{X_n\} = \{1, 1, 3, 7, 17, 41, 99, ...\}$.

From the recurrent dependencies for $Z_n$ and $Y_n$ we can determine that $Y_0 = 0$ and $Z_0 = 1$. It has become clear that $\{Y_n\} = \{0, 1, 2, 5, 12, 29, 70, ...\}$ are the famous Pell numbers ([2], **A000129)** whereas $\{Z_n\} = \{1, 3, 7, 17, 41, 99, ...\}$ and $\{X_n\} = \{1, 1, 3, 7, 17, 41, 99, ...\}$ are a sequence of numbers included in The On-Line Encyclopedia of Integer Sequences ([2], **A001333**). □

Let us examine the triangular lattices shown in Figure 10.

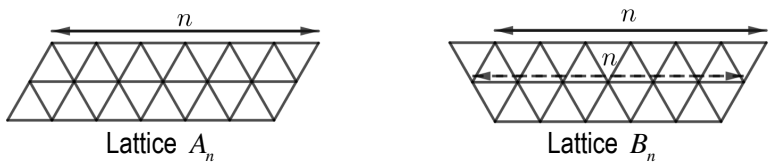

**Figure 10.** We consider the lattice $B_n$ to be formed from lattice $A_n$ by the moving of one of the small triangles - ($n$ is the midsegment of the trapezoid).

**Theorem 3.** The number of ways in which the lattices $A_n$ and $B_n$ can be covered with a combination of the two types of triangular tiles is given by the formulas:

$A_n = \frac{1}{4}\left[\left(1+\sqrt{2}\right)^{n+1} + \left(1-\sqrt{2}\right)^{n+1}\right] + \frac{1}{2}(-1)^n$, for $n \geq 1$,

$B_n = \frac{1}{4}\left[\left(1+\sqrt{2}\right)^{n+1} + \left(1-\sqrt{2}\right)^{n+1}\right] - \frac{1}{2}(-1)^n$, for $n \geq 1$,

$A_n = A_{n-1} + 3A_{n-2} + A_{n-3}$, for $n \geq 3$,

$B_n = B_{n-1} + 3B_{n-2} + B_{n-3}$, for $n \geq 3$,

$A_n = 2A_{n-1} + A_{n-2} + (-1)^n$, for $n \geq 2$,

$B_n = 2B_{n-1} + B_{n-2} + (-1)^n$, for $n \geq 2$.

**Proof.** First, we count directly $A_1 = 1$, $A_2 = 4$, $B_1 = 2$, $B_2 = 3$. This is illustrated in Figure 11 and Figure 12.

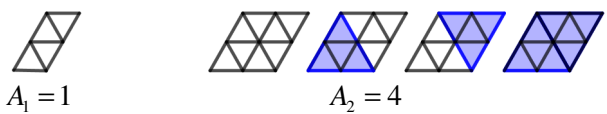

**Figure 11**

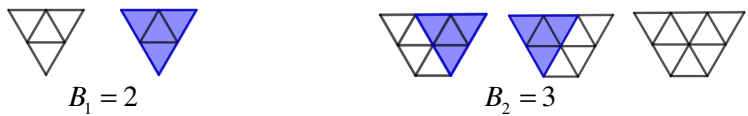

**Figure 12.**

Using the results from Theorem 1 and analyzing Figure 13, we conclude that for $n \geq 3$

$A_n = B_{n-1} + P_n$, $B_n = A_{n-1} + P_n$.

Thus
$$A_n = A_{n-2} + P_{n-1} + P_n = A_{n-2} + H_n,$$
$$B_n = B_{n-2} + P_{n-1} + P_n = B_{n-2} + H_n.$$

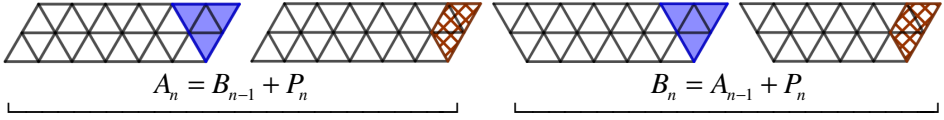

$$A_n = B_{n-1} + P_n \qquad B_n = A_{n-1} + P_n$$

**Figure 13**. The equations $A_n = B_{n-1} + P_n$, $B_n = A_{n-1} + P_n$ are illustrated.

Using recurrent equations we determine that $A_0 = 1$, $B_0 = 0$.
We find the generating functions and explicit formulas:
$$G_A = \sum_{n=0}^{\infty} A_n x^n = A_0 + A_1 x + \sum_{n=2}^{\infty} A_n x^n = 1 + x + \sum_{n=2}^{\infty} (A_{n-2} + H_n) x^n$$
$$= 1 + x + x^2 G_A + G_H - H_1 x - H_0 = x^2 G_A + \frac{1-x}{1-2x-x^2},$$
$$G_A = \frac{1}{(1+x)(1-2x-x^2)} = \frac{1}{1-x-3x^2-x^3}.$$

From the generating functions we get the recurrent equations:
$$A_n = A_{n-1} + 3A_{n-2} + A_{n-3},$$
from which follows the explicit formula for $n \geq 1$;
$$A_n = \frac{1}{4}\left[(1+\sqrt{2})^{n+1} + (1-\sqrt{2})^{n+1}\right] + \frac{1}{2}(-1)^n.$$
$$A_n = \frac{1}{4}H_{n+1} + \frac{1}{2}(-1)^n.$$
The recurrent equations for $n \geq 3$:
$$A_n = 2A_{n-1} + A_{n-2} + (-1)^n.$$
Similarly, we get the following results for the sequence $B_n$:
$$G_B = \frac{2x + x^2}{(1+x)(1-2x-x^2)} = \frac{2x+x^2}{1-x-3x^2-x^3},$$
$$B_n = B_{n-1} + 3B_{n-2} + B_{n-3},$$
$$B_n = \frac{1}{4}\left[(1+\sqrt{2})^{n+1} + (1-\sqrt{2})^{n+1}\right] - \frac{1}{2}(-1)^n,$$
$$B_n = \frac{1}{4}H_{n+1} - \frac{1}{2}(-1)^n, \quad B_n = 2B_{n-1} + B_{n-2} + (-1)^n.$$
It has become clear that
$$\{A_n\} = \{1, 4, 8, 21, 49, 120, 288, 697, 1681, 4060, \ldots\},$$
$$\{B_n\} = \{2, 3, 9, 20, 50, 119, 289, 696, 1682, 4059, \ldots\}.$$

Here is how the sequences mentioned above are presented in The On-Line Encyclopedia of Integer Sequences:

$\{A_n\}$ (виж **OEIS** [2], **A097076**),

0, 1, 1, 4, 8, 21, 49, 120, 288, 697, 1681, 4060, 9800, 23661, 57121, ...

and $\{B_n\}$ (виж **OEIS** [2], **A097075**),

1, 0, 2, 3, 9, 20, 50, 119, 289, 696, 1682, 4059, 9801, 23660, 57122, ... □

**Consequence.**

(1)   $A_{2n} = 6A_{2n-2} - A_{2n-4} - 2$  (**check out OEIS** [2], **A046090**),

1, 4, 21, 120, 697, 4060, 23661, 137904, ....

(2)   $A_{2n-1} = 6A_{2n-1} - A_{2n-3} + 2$  (**check out OEIS** [2], **A001108**),

0, 1, 8, 49, 288, 1681, 9800, 57121, ....

(3)   $B_{2n} = 6B_{2n-2} - B_{2n-4} + 2$  (**check out OEIS** [2], **A001652**),

0, 3, 20, 119, 696, 4059, 23660, 137903,...

(4)   $B_{2n-1} = 6B_{2n-1} - B_{2n-3} - 2$  (**check out OEIS** [2], **A055997**),

1, 2, 9, 50, 289, 1682, 9801, 57122,...

There are many dependencies among the tilings of the lattices $H_n$, $P_n$, $A_n$ and $B_n$. We are going to provide proof for some of them.

**Theorem 4.** For the lattices $H_n$, $P_n$, $A_n$ and $B_n$ the following equations are true for $n \geq 3$:

(1)   $A_n = B_n + (-1)^n$

(2)   $A_n = B_{n-1} + P_n$, $B_n = A_{n-1} + P_n$,

(3)   $A_{n-1} + A_{n-2} = B_{n-1} + B_{n-2} = P_n$,

(4)   $A_n = 2B_{n-1} + B_{n-2}$, $B_n = 2A_{n-1} + A_{n-2}$,

(5)   $A_n - A_{n-2} = H_n$, $B_n - B_{n-2} = H_n$,

(6)   $A_{2n} = A_0 + \sum_{i=1}^{n} H_{2i}$,  $A_{2n+1} = A_1 + \sum_{i=1}^{n} H_{2i+1}$

(7)   $B_{2n} = B_0 + \sum_{i=1}^{n} H_{2i}$,  $B_{2n+1} = B_1 + \sum_{i=1}^{n} H_{2i+1}$.

(8)   $P_n = P_{n-1} + A_{n-1} + B_{n-1}$.

**Proof.**

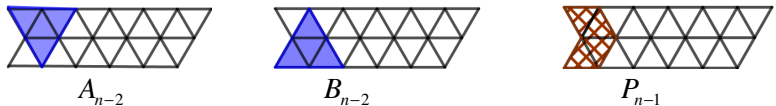

**Figure 14.** $P_n = A_{n-2} + B_{n-2} + P_{n-1}$.

For the equations (3) we are going to use (2) and the identity illustrated in Figure 14:

$P_n = A_{n-2} + (B_{n-2} + P_{n-1}) = A_{n-2} + A_{n-1}$,

$$P_n = A_{n-2} + B_{n-2} + P_{n-1} = (A_{n-2} + P_{n-1}) + B_{n-2} = B_{n-1} + B_{n-2}.$$

The result we got is the needed proof. For (4) we use (2), (3):
$$A_n = B_{n-1} + P_n = B_{n-1} + B_{n-1} + B_{n-2} = 2B_{n-1} + B_{n-2}.$$

For (5) we use the recurrent equations we have proven to be correct in Theorem 3.

We get the formulas (6) and (7) from the identities (5) by term-wise addition. For example, for the first formula we get the sum of the equations

$$A_2 - A_0 = H_2,\ A_4 - A_2 = H_4,\ A_6 - A_4 = H_6,\ \dots\ A_{2n} - A_{2n-2} = H_{2n}.\ \square$$

**Theorem 5.** For the lattices $H_n$, $P_n$, $A_n$ and $B_n$ the following equations are true for $n \geq 3$:

(1) $2P_{2n}^2 = A_{4n-1} = B_{4n-1} - 1$,

(2) $2P_{2n+1}^2 = B_{4n+1} = A_{4n+1} + 1$,

(3) $H_{2n}^2 = B_{4n-1} = A_{4n-1} + 1$,

(4) $H_{2n+1}^2 = A_{4n+1} = B_{4n+1} - 1$,

(5) $2P_{2n}P_{2n+1} = B_{4n} = A_{4n} - 1$,

(6) $2P_{2n+1}P_{2n+2} = A_{4n+2} = B_{4n+2} + 1$,

(7) $H_{2n}H_{2n+1} = A_{4n} = B_{4n} + 1$,

(8) $H_{2n+1}H_{2n+2} = B_{4n+2} = A_{4n+2} - 1$.

| $n$ | $H_n$ | $P_n$ | $A_n$ | | $B_n$ | |
|---|---|---|---|---|---|---|
| 1 | 1 | 1 | $1^2$ | 1 | $2.1^2$ | 2 |
| 2 | 3 | 2 | 2.1.2 | 4 | 1.3 | 3 |
| 3 | 7 | 5 | $2.2^2$ | 8 | $3^2$ | 9 |
| 4 | 17 | 12 | 3.7 | 21 | 2.2.5 | 20 |
| 5 | 41 | 29 | $7^2$ | 49 | $2.5^2$ | 50 |
| 6 | 99 | 70 | 2.5.12 | 120 | 7.17 | 119 |
| 7 | 239 | 169 | $2.12^2$ | 288 | $17^2$ | 289 |
| 8 | 577 | 408 | 17.41 | 697 | 2.12.29 | 696 |
| 9 | 1393 | 985 | $41^2$ | 1681 | $2.29^2$ | 1682 |
| 10 | 3363 | 2378 | 2.29.70 | 4060 | 41.99 | 4059 |
| 11 | 8119 | 5741 | $2.70^2$ | 9800 | $99^2$ | 9801 |
| 12 | 19601 | 13860 | 99.239 | 23661 | 2.70.169 | 23660 |

**Table 1.** The number of tilings for the lattices will be denoted exactly the same as the lattices themselves: $H_n$, $P_n$, $A_n$, $B_n$.

**Proof.**

(1) $\quad 2P_{2n}^2 = 2\left\{\dfrac{1}{2\sqrt{2}}\left[\left(1+\sqrt{2}\right)^{2n} - \left(1-\sqrt{2}\right)^{2n}\right]\right\}^2$

$$= 2\left\{\frac{1}{8}\left[\left(1+\sqrt{2}\right)^{4n}+\left(1-\sqrt{2}\right)^{4n}\right]-\frac{1}{4}(-1)^{2n}\right\} = A_{4n-1} = B_{4n-1}-1.$$

(2) $\quad 2P_{2n+1}^2 = 2\left\{\frac{1}{2\sqrt{2}}\left[\left(1+\sqrt{2}\right)^{2n+1}-\left(1-\sqrt{2}\right)^{2n+1}\right]\right\}^2$

$$= 2\left\{\frac{1}{8}\left[\left(1+\sqrt{2}\right)^{4n+2}+\left(1-\sqrt{2}\right)^{4n+2}\right]-\frac{1}{4}(-1)^{2n+1}\right\} = B_{4n+1} = A_{4n+1}+1.$$

(3) $\quad H_{2n}^2 = \left\{\frac{1}{2}\left[\left(1+\sqrt{2}\right)^{2n}+\left(1-\sqrt{2}\right)^{2n}\right]\right\}^2$

$$= \frac{1}{4}\left[\left(1+\sqrt{2}\right)^{4n}+\left(1-\sqrt{2}\right)^{4n}\right]+\frac{1}{2}(-1)^{2n} = B_{4n-1} = A_{4n-1}+1.$$

(4) $\quad H_{2n+1}^2 = \left\{\frac{1}{2}\left[\left(1+\sqrt{2}\right)^{2n+1}+\left(1-\sqrt{2}\right)^{2n+1}\right]\right\}^2$

$$= \frac{1}{4}\left[\left(1+\sqrt{2}\right)^{4n+2}+\left(1-\sqrt{2}\right)^{4n+2}\right]+\frac{1}{2}(-1)^{2n+1} = A_{4n+1} = B_{4n+1}-1.$$

(5) $\quad 2P_{2n}P_{2n+1} = 2\frac{1}{2\sqrt{2}}\left[\left(1+\sqrt{2}\right)^{2n}-\left(1-\sqrt{2}\right)^{2n}\right]\frac{1}{2\sqrt{2}}\left[\left(1+\sqrt{2}\right)^{2n+1}-\left(1-\sqrt{2}\right)^{2n+1}\right]$

$$= \frac{1}{4}\left[\left(1+\sqrt{2}\right)^{4n+1}+\left(1-\sqrt{2}\right)^{4n+1}\right]-\frac{1}{2} = B_{4n} = A_{4n}-1.$$

(6) $\quad 2P_{2n+1}P_{2n+2} = 2\frac{1}{2\sqrt{2}}\left[\left(1+\sqrt{2}\right)^{2n+1}-\left(1-\sqrt{2}\right)^{2n+1}\right]\frac{1}{2\sqrt{2}}\left[\left(1+\sqrt{2}\right)^{2n+2}-\left(1-\sqrt{2}\right)^{2n+2}\right]$

$$= \frac{1}{4}\left[\left(1+\sqrt{2}\right)^{4n+2}+\left(1-\sqrt{2}\right)^{4n+2}\right]-\frac{1}{2} = A_{4n+2} = B_{4n+2}+1.$$

(7) $\quad H_{2n}H_{2n+1} = \frac{1}{2}\left[\left(1+\sqrt{2}\right)^{2n}+\left(1-\sqrt{2}\right)^{2n}\right]\frac{1}{2}\left[\left(1+\sqrt{2}\right)^{2n+1}+\left(1-\sqrt{2}\right)^{2n+1}\right]$

$$= \frac{1}{4}\left[\left(1+\sqrt{2}\right)^{4n+1}+\left(1-\sqrt{2}\right)^{4n+1}\right]+\frac{1}{2} = A_{4n} = B_{4n}+1.$$

(8) $\quad H_{2n+1}H_{2n+2} = \frac{1}{2}\left[\left(1+\sqrt{2}\right)^{2n+1}+\left(1-\sqrt{2}\right)^{2n+1}\right]\frac{1}{2}\left[\left(1+\sqrt{2}\right)^{2n+2}+\left(1-\sqrt{2}\right)^{2n+2}\right]$

$$= \frac{1}{4}\left[\left(1+\sqrt{2}\right)^{4n+3}+\left(1-\sqrt{2}\right)^{4n+3}\right]-\frac{1}{2} = B_{4n+2} = A_{4n+2}-1. \quad \square$$

**Задача 1.** An angular lattice consists of two triangular lattices $X_n$, as shown in Figure 15. The dimensions of the given grid are $2\times n\times n$, $n\geq 1$. Find the number of different ways in which the lattice, depicted in Figure 15, can be covered using a combination of tiles of the two types of triangles.

**Solution.** Let us denote the number of different ways in which the grid in Figure 15 can be covered with triangles with $C_n$. First, we directly count $C_1 = 1$, $C_2 = 3$.

After dividing the lattice into two lattices, as shown in Figure 17, we see that part of the figure can be covered in $X_n^2$ ways and the other part of the figure can be covered in $2Y_{n-1}Z_{n-1}$ ways.

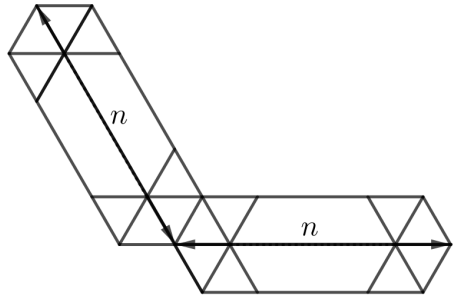

**Figure 15**

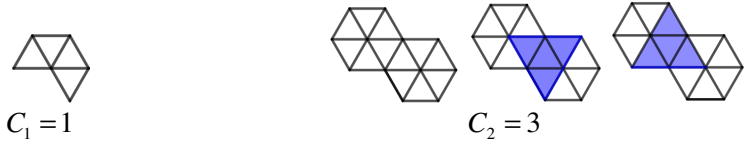

$C_1 = 1$          $C_2 = 3$

**Figure 16.**

Let us consider the triangular lattices depicted in Figure 17:
$C_1 = 1$, $C_2 = X_2^2 + 2Y_1Z_1 = 1 + 2.1.1 = 3$,

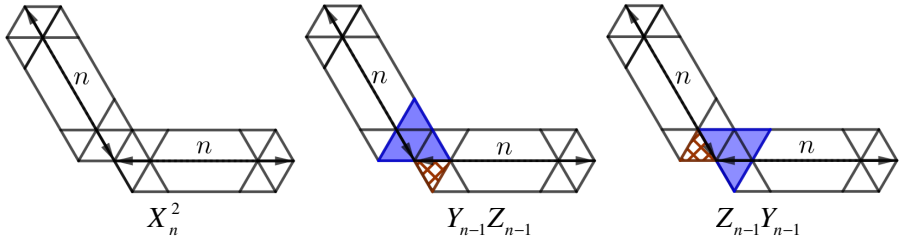

$X_n^2$          $Y_{n-1}Z_{n-1}$          $Z_{n-1}Y_{n-1}$

**Figure 17.** The equation $C_n = X_n^2 + 2Y_{n-1}Z_{n-1}$ is illustrated.

$C_n = X_n^2 + 2Y_{n-1}Z_{n-1} = H_{n-1}^2 + 2P_{n-1}H_{n.1} = H_{n-1}^2 + P_{2n-2}$ (according to $2P_nH_n = P_{2n}$).
According to Theorem 4(2) $A_n = B_{n-1} + P_n$, $B_n = A_{n-1} + P_n$ and to Theorem 5(3) $H_{2n}^2 = B_{4n-1}$, 5(4) $H_{2n+1}^2 = A_{4n+1}$:
$C_{2n-1} = B_{4n-1} + P_{4n-2} = A_{4n}$,
$C_{2n} = A_{4n-3} + P_{4n-2} = B_{4n-2}$.
Now we reach the following dependencies:
$C_{n+1} = H_n^2 + P_{2n} =$

$$= \frac{1}{4}\left(1+\sqrt{2}\right)^{2n+1} + \frac{1}{4}\left(1-\sqrt{2}\right)^{2n+1} + \frac{1}{2}(-1)^n = \frac{1}{2}H_{2n+1} + \frac{1}{2}(-1)^n.$$

Here is how the lattice $C_n$ is presented in **OEIS, (A084159,** [2]):

1, 3, 21, 119, 697, 4059, 23661, 137903,...

The solutions here are also a sequence with recurrence relations:

$C_{2n+1} = 6C_{2n} - C_{2n-1} + 4$, $n \geq 1$,

$C_{2n} = 6C_{2n-1} - C_{2n-2} - 4$, $n \geq 2$. □

**Problem 2.** An angular grid consists of two triangular grids $Y_n$, as shown in Figure 18. The dimensions of the given grid are $2 \times n \times n$, $n \geq 1$. Find the number of different ways in which the lattice, depicted in Figure 18, can be covered using a combination of tiles of the two types of triangles.

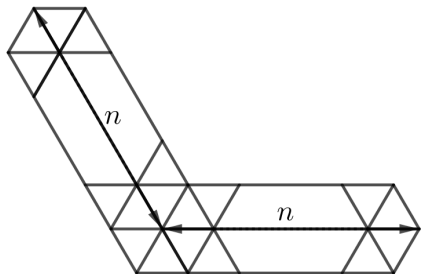

**Figure 18**

**Solution.** Let us denote the number of different ways in which the grid in Figure 18 can be covered with triangles with $D_n$. First, we directly count $D_1 = 1$, $D_2 = 6$.

After dividing the lattice into two lattices, as shown in Figure 20, we see that part of

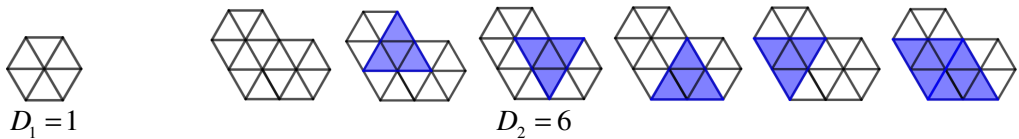

$D_1 = 1$ $\qquad\qquad\qquad D_2 = 6$

**Figure 19**

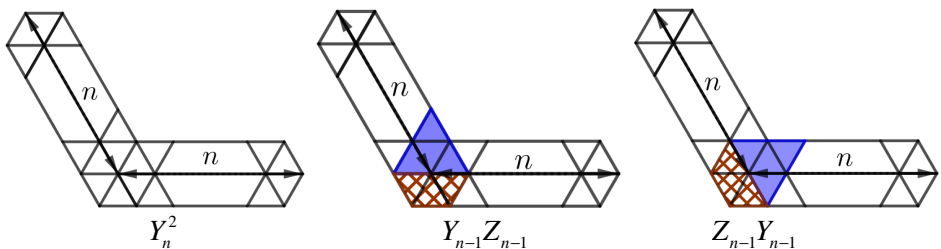

$Y_n^2 \qquad\qquad Y_{n-1}Z_{n-1} \qquad\qquad Z_{n-1}Y_{n-1}$

**Figure 20.** The equation $D_n = Y_n^2 + 2Y_{n-1}Z_{n-1}$ is illustrated.

the figure can be covered in $Y_n^2$ ways and the other part of the figure can be covered in $2Y_{n-1}Z_{n-1}$ ways.

Let us consider the triangular lattices depicted in Figure 20:
$D_1 = X_1 = 1$, $D_2 = P_2^2 + 2P_1 H_1 = 4 + 2.1.1 = 6$,
$D_n = Y_n^2 + 2Y_{n-1}Z_{n-1} = P_n^2 + 2P_{n-1}H_{n-1} = P_n^2 + P_{2n-2}$

$$= \left\{ \frac{1}{2\sqrt{2}} \left[ (1+\sqrt{2})^n - (1-\sqrt{2})^n \right] \right\}^2 + \frac{1}{2\sqrt{2}} \left[ (1+\sqrt{2})^{2n-2} - (1-\sqrt{2})^{2n-2} \right]$$

$$= \frac{3+4\sqrt{2}}{8}(1+\sqrt{2})^{2n-2} + \frac{3-4\sqrt{2}}{8}(1-\sqrt{2})^{2n-2} - \frac{1}{4}(-1)^n.$$

Thus we get the formula
$$D_n = \frac{3+4\sqrt{2}}{8}(1+\sqrt{2})^{2n-2} + \frac{3-4\sqrt{2}}{8}(1-\sqrt{2})^{2n-2} - \frac{1}{4}(-1)^n.$$

The obtained formula leads to the conclusion that the sequence $D_n$ has the characteristic equation $(x^2 - 6x + 1)(x+1) = 0$.

The solutions here are also a sequence with recurrence relations:
$D_{2n} = 6D_{2n-1} - D_{2n-2} - 2$, $n \geq 2$,
$D_{2n+1} = 6D_{2n} - D_{2n-1} + 2$, $n \geq 1$.

**Note.** The sequence $\{D_n\}$ is not listed in The On-Line Encyclopedia of Integer Sequences: 1, 6, 37, 214, 1249, 7278, 42421, ...

**Let us consider the grids $H_n$ and $P_n$ again**, which are depicted in Figure 3, denoting respectively with $h_n$ and $p_n$ the total number of tiles (small and large) needed for the realization of all the tilings, with $\varphi_n$ and $\theta_n$ the number of only the small tiles, and with $q_n$ and $r_n$ the number of only the large tiles.

| n | $H_n$ 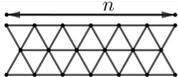 Lattice $H_n$ | | | | $P_n$ 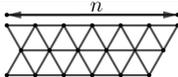 Lattice $P_n$ | | | |
|---|---|---|---|---|---|---|---|---|
| | $H_n$ | $h_n$ | $\varphi_n$ | $q_n$ | $P_n$ | $p_n$ | $\theta_n$ | $r_n$ |
| 1 | 1 | 2 | 2 | 0 | 1 | 1 | 1 | 0 |
| 2 | 3 | 12 | 10 | 2 | 2 | 7 | 6 | 1 |
| 3 | 7 | 46 | 38 | 8 | 5 | 30 | 25 | 5 |
| 4 | 17 | 154 | 126 | 28 | 12 | 102 | 84 | 18 |
| 5 | 41 | 474 | 390 | 88 | 29 | 319 | 257 | 58 |
| 6 | 99 | 1392 | 1138 | 262 | 70 | 945 | 762 | 175 |
| 7 | 239 | 3958 | 3226 | 752 | 169 | 2704 | 2177 | 507 |
| 8 | 577 | 10692 | 8942 | 2104 | 408 | 7548 | 6072 | 1428 |

**Table 2.** The values for the considered series are shown.

**Theorem 6. a)** The number of tiles (small and large) needed for the tilings of lattices $H_n$ and $P_n$, shown in Figure 3, is expressed with the following formulas for $n \geq 1$:

$$h_n = \frac{5}{2}nH_n + \frac{3}{2}P_n - 2H_n, \quad p_n = \frac{5}{2}nP_n - \frac{3}{2}P_n.$$

**б)** The number of large tiles needed for the tilings of lattices $H_n$ and $P_n$, shown in Figure 3, is expressed with the following formulas for $n \geq 1$:

$$q_n = \frac{1}{2}nH_n - \frac{1}{2}P_n, \quad r_n = \frac{1}{2}nP_n - \frac{1}{2}P_n.$$

**в)** The number of small tiles needed for the tilings of lattices $H_n$ and $P_n$, shown in Figure 3, is expressed with the following formulas for $n \geq 1$:

$$\varphi_n = 2(n-1)H_n + 2P_n, \quad \theta_n = 2nP_n - P_n = (2n-1)P_n.$$

**Proof.** We will derive recurrence dependencies regarding the number of tiles for both lattices by comparing the number of necessary tiles for all tilings of a strip of size $n$ and a strip of size $n-1$.

The following equations are illustrated in Figure 21:
$h_n = h_{n-1} + 4H_{n-1} + 2p_{n-1} + 4P_{n-1}$,
$q_n = q_{n-1} + 2r_{n-1} + 2P_{n-1}$,
$\varphi_n = \varphi_{n-1} + 4H_{n-1} + 2\theta_{n-1} + 2P_{n-1}$.

On the first diagram in Figure 21, the case is shown where none of the four added small triangles on the strip of size $n$, with respect to a strip of size $n-1$, has a large tile placed on it. In this case, to the number $h_{n-1}$ of the necessary tiles for all tilings on the strip of size $n$ four small tiles are added for each tiling on the strip of size $n-1$, adding $4H_{n-1}$ tiles in total.

The next two diagrams in Figure 21 illustrate the cases where one of the rightmost small triangles is covered by a large tile. Then, to the number $p_{n-1}$ of necessary tiles for all tilings on the lattice of size $n-1$ one small and one large tile are added for each tiling on the strip of size $n-1$, adding $2P_{n-1}$ tiles in total. Ultimately for $n \geq 2$

$h_n = h_{n-1} + 4H_{n-1} + 2(p_{n-1} + 2P_{n-1})$.

Similar reasoning applies to the remaining equations in Figure 21.

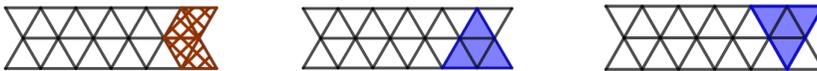

| | | | | | |
|---|---|---|---|---|---|
| $h_n =$ | $h_{n-1} + 4H_{n-1}$ | $+$ | $p_{n-1} + 2P_{n-1}$ | $+$ | $p_{n-1} + 2P_{n-1}$ |
| $q_n =$ | $q_{n-1}$ | $+$ | $r_{n-1} + P_{n-1}$ | $+$ | $r_{n-1} + P_{n-1}$ |
| $\varphi_n =$ | $\varphi_{n-1} + 4H_{n-1}$ | $+$ | $\theta_{n-1} + P_{n-1}$ | $+$ | $\theta_{n-1} + P_{n-1}$ |

**Figure 21.** The equations $h_n$, $q_n$, $\varphi_n$ are illustrated on the diagrams.

The following equations for $n \geq 2$ are illustrated in Figure 22:

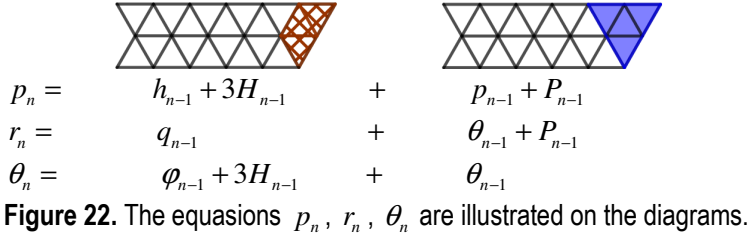

$$p_n = \quad h_{n-1}+3H_{n-1} \quad + \quad p_{n-1}+P_{n-1}$$
$$r_n = \quad q_{n-1} \quad + \quad \theta_{n-1}+P_{n-1}$$
$$\theta_n = \quad \varphi_{n-1}+3H_{n-1} \quad + \quad \theta_{n-1}$$

**Figure 22.** The equasions $p_n$, $r_n$, $\theta_n$ are illustrated on the diagrams.

$p_n = h_{n-1}+3H_{n-1}+p_{n-1}+P_{n-1}$,
$r_n = q_{n-1}+r_{n-1}+P_{n-1}$,
$\theta_n = \varphi_{n-1}+3H_{n-1}+\theta_{n-1}$.

Using recurrence dependencies illustrated in Figures 21 and 22, we determine the number of tilings for $n \geq 3$:

$h_n = 2h_{n-1}+h_{n-2}+10P_{n-1}$,
$p_n = 2p_{n-1}+p_{n-2}+5H_{n-1}$,
$q_n = 2q_{n-1}+q_{n-2}+2P_{n-1}$,
$r_n = 2r_{n-1}+r_{n-2}+H_{n-1}$,
$\varphi_n = 2\varphi_{n-1}+\varphi_{n-2}+8P_{n-1}$,
$\theta_n = 2\theta_{n-1}+\theta_{n-2}+4H_{n-1}$.

For the total number of tiles:
**a)** We form a system of two recurrent equations

$$\begin{vmatrix} h_n = h_{n-1}+2p_{n-1}+4H_{n-1}+4P_{n-1} \\ p_n = p_{n-1}+h_{n-1}+P_{n-1}+3H_{n-1} \end{vmatrix}.$$

Finally, through recurrent relations we get:
$2p_{n-1} = h_n - h_{n-1} - 4P_n$ и $h_{n-1} = p_n - p_{n-1} - P_{n-1} - 3H_{n-1}$

Then we substitute in the system and we get two non-homogeneous equations
$h_n - 2h_{n-1} - h_{n-2} = 10P_{n-1}$, $n \geq 2$;
$p_n - 2p_{n-1} - p_{n-2} = 5H_{n-1}$, $n \geq 2$.

They have the same characteristic equation
$x^2 - 2x - 1 = 0$.

With direct verification and recurrent dependencies we obtain:
$h_1 = 2$, $h_2 = 12$, $h_3 = 46$, $h_3 = 154$;
$p_1 = 1$, $p_2 = 7$, $p_3 = 30$, $p_4 = 102$.

We derive the explicit formula for the number $h_n$ of all tiles for all different tilings of lattice $H_n$ in Figure 21 as follows: Given the expression for $P_{n-1}$ on the right side of the non-homogeneous recurrence equation $h_n - 2h_{n-1} - h_{n-2} = 10P_{n-1}$, the wanted formula for $n \geq 1$ is of the following form:
$$h_n = c_1\left(1+\sqrt{2}\right)^n + c_2\left(1-\sqrt{2}\right)^n + c_3 n\left(1+\sqrt{2}\right)^n + c_4 n\left(1-\sqrt{2}\right)^n.$$

The constants are found through the system obtained by substituting $n$ with the numbers 1, 2, 3, 4 sequentially:

Ultimately

$$h_n = \frac{3}{4\sqrt{2}}\left[\left(1+\sqrt{2}\right)^n - \left(1-\sqrt{2}\right)^n\right] - \left[\left(1+\sqrt{2}\right)^n - \left(1-\sqrt{2}\right)^n\right] + \frac{5}{4}n\left[\left(1+\sqrt{2}\right)^n + \left(1-\sqrt{2}\right)^n\right].$$

By transforming the explicit formula for the number of all tiles for all different tilings of lattice $H_n$ we obtain an impressive and also convenient formula for $h_n$:

$$h_n = \frac{5}{2}nH_n + \frac{3}{2}P_n - 2H_n = \frac{5n+1}{2}H_n - \frac{3}{2}H_{n-1}.$$

Similarly, we derive the formulas for the number of all tiles for all tilings of lattice $P_n$.

$$2p_{n-1} = h_n - h_{n-1} - 4P_n,$$

$$2p_n = \frac{5}{2}(n+1)H_{n+1} + \frac{3}{2}P_{n+1} - 2H_{n+1} - \left(\frac{5}{2}nH_n + \frac{3}{2}P_n - 2H_n\right) - 4P_{n+1},$$

$$p_n = \frac{5}{2}nP_n - \frac{3}{2}P_n = \frac{5n-3}{2}P_n = \frac{5n-3}{2}\left(H_n - H_{n-1}\right).$$

**б)** For the number of large tiles $q_n$ for all tilings of lattices $H_n$ we proved the non-homogeneous equation

$q_n - 2q_{n-1} - q_{n-2} = 2P_n$, $n \geq 2$, and we can proceed similarly - the sought formula has the following form for $n \geq 1$:

$$q_n = c_1\left(1+\sqrt{2}\right)^n + c_2\left(1-\sqrt{2}\right)^n + c_3 n\left(1+\sqrt{2}\right)^n + c_4 n\left(1-\sqrt{2}\right)^n.$$

Here, too, the constants are found through the system obtained by substituting $n$ sequentially with the numbers 1, 2, 3, 4:

Ultimately

$$q_n = -\frac{1}{4\sqrt{2}}\left[\left(1+\sqrt{2}\right)^n - \left(1-\sqrt{2}\right)^n\right] + \frac{1}{4}n\left[\left(1+\sqrt{2}\right)^n + \left(1-\sqrt{2}\right)^n\right].$$

From the proven formulas for $h_n$ and $q_n$ it has become evident that

$$q_n = \frac{1}{2}nH_n - \frac{1}{2}P_n.$$

For the number of large tiles for all tilings of lattices $P_n$ - the sought formula has the following form for $n \geq 1$:

$$q_n = q_{n-1} + 2r_{n-1} + 2P_{n-1},$$

$$2r_n = \frac{1}{2}(n+1)H_{n+1} - \frac{1}{2}P_{n+1} - \left(\frac{1}{2}nH_n - \frac{1}{2}P_n\right) - 2P_n.$$

Now it is easy to see that

$$r_n = \frac{1}{2}nP_n - \frac{1}{2}P_n,$$

$$q_n = \frac{1}{2}nH_n - \frac{1}{2}P_n.$$

It is evident that $\{r_n\}$ are **OEIS** ($[2]$, **A364553**) – a sequence of numbers
0, 1, 5, 18, 58, 175, 507, 1428, 3940, ...

It is evident that $\left\{\dfrac{q_n}{2}\right\}$ are **OEIS** ($[2]$, **A006645**) – a sequence of numbers
0, 1, 4, 14, 44, 131, 376, 1052, 2888, 7813, ...

в) For the number of small tiles for all tilings of lattices $H_n$ and $P_n$, which are shown in Figure 21, we can use the already proven formulas: $\varphi_n = h_n - q_n$, $\theta_n = p_n - r_n$.

From here we have
$\varphi_n = 2nH_n + 2P_n - 2H_n = 2(n-1)H_n + 2P_n$,
$\theta_n = 2nP_n - P_n = (2n-1)P_n$. $\square$

**Theorem 7.** Let us consider lattices $H_n$ and $P_n$, which are shown in Figure 3, denoting respectively with $h_n$ and $p_n$ the total number of tiles (small and large) for the realization of all tilings, with $\varphi_n$ and $\theta_n$ the number of only the small tiles. Then
$$\lim_{n\to\infty}\frac{\varphi_n}{h_n}=\frac{4}{5},\ \lim_{n\to\infty}\frac{\theta_n}{p_n}=\frac{4}{5}.$$

**Proof.** We will use the equations:
$h_n = \dfrac{5n}{2}H_n + \dfrac{1}{2}H_n + \dfrac{3}{2}P_n$,
$\varphi_n = 2nH_n + 2P_n$.

Since $h_n \sim \dfrac{5n}{2}H_n$ and $\varphi_n \sim 2nH_n$, then $\lim_{n\to\infty}\dfrac{\varphi_n}{h_n}=\dfrac{4}{5}$.

Respectively, from these equations we get
$p_n = \dfrac{5n+2}{2}P_n$,
$\theta_n = (2n+1)P_n$.

That means $p_n \sim \dfrac{5n}{2}P_n$ and $\theta_n \sim 2nP_n$, so $\lim_{n\to\infty}\dfrac{\theta_n}{p_n}=\dfrac{4}{5}$. $\square$

Similarly, we derive the formulas $\lim_{n\to\infty}\dfrac{q_n}{h_n}=\dfrac{1}{5}$, $\lim_{n\to\infty}\dfrac{r_n}{p_n}=\dfrac{1}{5}$.

Let us consider the lattices $A_n$ and $B_n$, depicted in Figure 10. Similarly to Theorem 6 and Theorem 7, we determine the total number of tiles (small and large) required for all tilings. Then, the asymptotics are also $\dfrac{4}{5}, \dfrac{1}{5}$.

Valcho Milchev

Kardzhali, Bulgaria

e-mail:
milchev.vi@gmail.com

milchev_v@abv.bg